\documentclass[12pt]{article}
\usepackage{amsmath,amsfonts, amssymb,verbatim}
\newtheorem{pkt}{}[section]
\newcommand{\bpk}{\begin{pkt}\rm }
\newcommand{\epk}{\end{pkt}}

\newcommand{\QQ}{{^*\mathbb{Q}}}

\newcommand{\inv}{^{-1}}

\newcommand{\hh}{\mathfrak{h}}
\newcommand{\M}{{\bf M}}

\newcommand{\Pp}{\mathrm{P}}
\newcommand{\Qq}{\mathrm{Q}}

\newcommand{\R}{{\mathbb R}}
\newcommand{\Q}{{\mathbb Q}}
\newcommand{\Z}{{\mathbb Z}}
\newcommand{\N}{{\mathbb N}}
\newcommand{\C}{{\mathbb C}}

\newcommand{\be}{\begin{equation}}
\newcommand{\ee}{\end{equation}}
\newcommand{\kk}{\mathrm{k}}

\newcommand{\Ss}{\mathbb{S}}

\newcommand{\K}{{\rm K}}

\newcommand{\NN}{\mathfrak{n}}
\newcommand{\MM}{\mathcal{M}}

\newcommand{\qed}{$\Box$ \\ \\}
\newcommand{\F}{\mathrm{F}}

\newcommand{\pp}{\mathfrak{p}}
\newcommand{\lfr}{\mathfrak{l}}

\newcommand{\ZZ}{{^*\mathbb{Z}}}

\newcommand{\lm}{\mathsf{lm}}

\newcommand{\e}{\mathrm{e}}

\title{ Approximation of structures: local and global}
\author{Boris Zilber}
\begin{document}
	\maketitle
		\abstract{We provide a mathematically rigorous definition of local approximation and consider it along with previously used global approximation.  We demonstrate the applicability of the new construction to some interesting classes of structures. 
		In particular, we prove that any compact simple Lie group is locally approximated by finite groups. The definition and main examples are motivated by physics but the techniques are of model theory. Namely, we introduce the ultraproduct of emerging metric structures, which generalises the ultraproduct in metric model theory. }
	\section{Introduction}
	\bpk  This work along with some earlier ones cited below is motivated by attempts to understand the procedures of approximation in physics, in particular by the question of how infinite and continuous emerges from  finite and discrete. It is an age-old question debated by Greek philosophers and reformulated by Hilbert  in his  problem 6 where he proposed  developing ``mathematically the limiting process ... which leads from the atomistic view to the laws of motion of continua'', see \cite{SR}.
	\epk
	\bpk 
	 We continue the study of structural approximation (of continuous structures by finite structures)  initiated in \cite{Perfect}. In  \cite{ZlbKL} and \cite{ZDvN} we realised that along with its initial {\em global} version there might be  another useful version which we identified in a few specific cases  and called {\em local} approximation.
	
	Here we furnish a mathematically rigorous definition of local approximation and show how it is applicable to some interesting classes of structures. This includes the approximation of the Lorentz-invariant Minkowski spacetime of \cite{ZlbKL} by finite  lattices  as well as the continuous model theory approximation of quantum mechanics from \cite{ZDvN}.

	 Also, addressing a physics motivated question posed in \cite{Perfect}, we show that local approximation of compact simple Lie groups by finite groups does exist, while 	 
	 \cite{Nik} proved that a global approximation is not possible.

	 Our definition of local approximation reduces to the definition of the {\em local ultraproduct}, the ultraproduct in the class of {\em  emerging metric} structures. This generalises the ultraproduct of metric structures. The construction is based on the idea that in some classes of  algebraic  structures such as groups or fields, in sufficiently large  finite structures there is a natural metric ``near the origin'' but which is  not necessarily well-defined globally. Under certain assumptions, in the first-order ultraproduct, the near-the-origin region will form an algebraic structure with a pseudo-metric, which then is reduced into a structure with a complete metric.

	 \medspace   
	 
	 Our context includes both first order topological structures as in \cite{Perfect} and \cite{ZCL} and metric  structures as in \cite{BBHU} as well as continuous logic structures as \cite{ZDvN}.  This becomes possible in the setting of {\em general structures} developed by Keisler, \cite{KeislerGL}, which allows any predicates, not necessarily continuous, with values in the real interval $[0,1].$ However, the role of the distance predicate is slightly different. In fact, our emerging metric is a generalisation of the  pre-metric expansion of Keisler.   
	
	\epk
	\bpk Applications of the new notion of approximation  in the paper demonstrate that  important structures of mathematical physics can be locally approximated by finite structures. In particular, in section~\ref{s4} we prove that {\em all compact simple Lie groups can be locally approximated by finite groups}. These groups are the backbone of the Standard Model of particle physics and also  of  string theory. However, it is worth remarking that from a calculational point of view the fact that a metric structure can be approximated  by finite emerging metric structures probably does not add much. The main outcome of the local approximation theory is that  the calculations are consistent with the assumption of finiteness of the model.  Global approximation, on the other hand, may be of importance when models of closed systems like a whole universe are considered, see, e.g. \cite{ZlbKL} and \cite{ZDvN} for speculations in this direction.
	
	In conclusion, we would like to formulate a relevant open problem. 
	
	{\bf The Virasoro group } is  the universal central extension of the infinite-dimensional Lie group of orientation preserving diffeomorphisms of the circle $\mathrm{Diff}^+(\Ss^1)$  by $\R,$  see \cite{Kir}. {\em Is the Virasoro group or just the group $\mathrm{Diff}^+(\Ss^1)$ locally approximable by finite groups? }
	
	Note that the Virasoro group can be considered a generalisation of the Heisenberg group, which is locally approximable, see \ref{Heis}.  
	\epk
	\bpk It is interesting to compare, as suggested by a referee, our notions of approximation with the definition studied in \cite{GGH} by  Glebsky, Gordon and Henson. Their notion of approximation by finite algebras is motivated by applications such as computer arithmetic, which leads to a technically  different definition. One of the main results of \cite{GGH} is the theorem stating that no locally compact field can be approximated by finite rings in their sense. In contrast, we establish, see \ref{3.4}, that the field of reals can be locally approximated by finite fields,  while the global approximation by finite fields is possible for the compactified field of complex numbers $\bar{\C}$ only.       
	\epk
	\bpk I am grateful to Arkady Bollotin, Simon Saunders, and Tim Palmer for fruitful discussions that helped motivate  the paper. 
	\epk
\section{Main definitions}
 \bpk\label{assm} {\bf Topological structures.}
 
 Let $\M$ be a  structure in a relational first-order language $L.$  
 The $L$-topology on $M^n,$ for all $n,$ is given by declaring a  subset $P\subseteq M^n$ closed if it is defined by positive $L$-formulas, that is by formulas built from basic relations using $\&,$ $\vee,$ $\exists$ and $\forall$ only.
 
 An $L$-structure $\M=(M;L)$ considered with $L$-topology will be referred to as a {\bf topological structure}.

 We say $\M$ is {\bf compact} (often  {\bf quasi-compact} in literature) if  every point in $M$ is closed and for 
 any filter of closed subsets of $M^n$ the intersection is nonempty.\\
 
 \epk

\bpk \label{Kies} {\bf The formalism of general structures} (following \cite{KeislerGL}). The space of truth values will be $[0, 1],$ with $0$ representing truth. A vocabulary  consists of a set of predicate symbols $P$ of finite arity,
a set of function symbols $f$ of finite arity, and a set of constant symbols $c.$ A general
 structure $\M$ consists of a vocabulary $V$, a non-empty universe set $M$,
an element $c^M\in M$ for each constant symbol $c$, a mapping $P^M : M^n \to [0, 1]$ for
each predicate symbol $P$ of arity $n$, and a mapping $f^M : M^n \to M$ for each function symbol $f$ of arity $n.$ A general structure determines a vocabulary, but does
not determine a metric signature.
The formulas are as in \cite{BBHU}, with the connectives  being all continuous functions from finite powers of $[0, 1]$ into $[0, 1],$ and the quantifiers $\sup x,$  $\inf x.$ The truth
value of a formula $\varphi(\bar{x})$ at a tuple  $\bar{a}\in M^n$ in a general structure $\M$ is an element
of $ [0, 1]$ denoted by $\varphi^\M(\bar{x}).$ It is defined in the usual way by induction on the
complexity of formulas. 

The case of interest to us is when the vocabulary $V=L\cup \{ \mathsf{d}  \},$ where $L$ contains names of $\{ 0,1 \}$ predicates $P$ only and  $\mathsf{d}: M^2\to [0,1]$ is the ``distance'' predicate (not necessarily satisfying the assumptions defining a metric).

We are going to apply these definitions in a many-sorted version, assuming that our sorts $S_m,$ $m\in \N ,$ are nested: $S_m\subseteq S_{m+1}$ (see \ref{L+}). The diameters of the sorts can grow, which is not consistent with the distance predicate taking values in $[0,1]$. This is easily remedied by allowing a separate distance predicate $\mathsf{d}_m$ for each sort $S_m$ and then postulating the scaling relation between $\mathsf{d}_m.$
\epk
\bpk\label{global} {\bf Global approximation of first-order topological structures.}  

Let $\M_{i\in I}$ be a family of topological  $L$-structures and $\M$ a topological $L$-structure of cardinality continuum.

Let $\mathcal{D}$ be a non-principal ultrafilter  on $I$ and let  $${\prod}_I\M_i/\mathcal{D} $$
be the first-order ultraproduct of $L$-structures. Let 
\be\label{M*} {^*\M}\succeq  {\prod}_I\M_i/\mathcal{D}\ee 
be a $2^{\aleph_0}$-saturated elementary extension of the ultraproduct.

Then, 
 $\M$ is  {\bf globally approximated} by the family $\M_{i\in I}$ along the ultrafilter $\mathcal{D}$ on $I$ if, for every ${^*\M}$ as in (\ref{M*}),
there is a surjective $L$-homomorphism
\be\label{lm2}\lm^\mathrm{glob}: {^*\M} \twoheadrightarrow \M.\ee

We often assume that  all  members of the family $\M_{i\in I}$ are finite models of a theory $T,$ in which case we say $\M$ is globally approximated by finite models of $T.$

We also write the fact of global approximation  in the form \be\label{lm}\lm^\mathrm{glob} \M_\mathbf{i}=\M\ee
which  omits  a reference to $I$ and the ultrafilter $\mathcal{D}.$

\epk
\bpk\label{sh} {\bf Strong homomorphisms.}  Let $\mathbf{N}$ and $\M$ be $L$-structures. 
 
Given $A\subset N,$ a partial homomorphism $h_A: \mathbf{N}\to \M$ is a map
$h_A:A\to \M$ 
defined on $A$ such that for every $a\in A^k$  and every  $L$-predicate $S(x)$  $$\mathbf{N}\vDash  S(a) \Rightarrow \M\vDash S(a^h),$$ where $a^h:=h_A(a).$
  
A {\em partial strong homomorphism} $h_A: \mathbf{N}\to \M$ is a map
 $h_A:A\to \M$ 
defined on $A$ such that for every $a\in A^k$  and every positive $L$-formula $S(x)$  $$\mathbf{N}\vDash  S(a) \Rightarrow \M\vDash S(a^h).$$   
This includes the case $A=\emptyset$, in which instance $S$ is just the sentence   and the definition requires $$\mathbf{N}\vDash S \Rightarrow \M\vDash S.$$
\epk
\bpk \label{PropGlob1} 
{\bf Proposition.} (
  Cf. \cite{Perfect}, Prop.3.4). {\em  If $\M$ is compact then the existence of ${^*\M}$ and  $\lm^\mathrm{glob}$ such that {\rm (\ref{lm2})}
   is equivalent to the assumption that for every positive
  $L$--sentence $\sigma$ 
  \be\label{LosPos}  {\prod}_I\M_i/\mathcal{D} \vDash \sigma \Rightarrow \M \vDash \sigma\ee
  (the positive-Lo\'s property). 

Moreover, {\rm (\ref{lm2})} implies that for any $\mathbf{N}\succeq {^*\M}$ there is a surjective homomorphism  
	 $$\lm: \mathbf{N}\twoheadrightarrow \M$$
 extending  $\lm^\mathrm{glob}.$ }

{\bf Proof.} (\ref{lm2}) implies (\ref{LosPos}) by the standard Lo\'s theorem. We need to show the inverse.
	
 First we note that (\ref{LosPos}) implies that there exists the partial strong homomorphism $h_\emptyset: \mathbf{N}\to \M.$ 

The proof proceeds with two claims.

Claim 1. Suppose for some $A\subset N$ there is a partial strong homomorphism
$h_A:A\to \M,$ and $b\in N.$ Then $h_A$ can be extended to a partial strong homomorphism $h_{Ab}:Ab\to \M.$

Proof of Claim. Let  $p^+_b(y)$ be a positive type of $b$ over $A$ in $\mathbf{N}.$ That is   $p^+_b(y)$ 
consists of positive formulas of the form $\exists z\, S(a,y,z),$ such that $\mathbf{N}\vDash \exists z\, S(a,b,z).$ 

Applying $h_A$ we get
  $$\M\vDash  \exists yz S(a^h,y,z),\mbox{ for all }\exists z S(a,y,z)\in p^+_b(y) $$ 
It follows that the family $\mathcal{P}^+_b(\M)$ of closed sets, each of the form $$\{ c\in M: \exists z S(a^h,c,z)\},\mbox{ for some }\exists z S(a,y,z)\in p^+_b(y)$$ 
is  
a filter. By quasi-compactness of $\M$ there is a point $b'$ in the intersection $\bigcap\mathcal{P}^+_b(\M).$  Clearly, letting $h_{Ab}: b\mapsto b':=b^h,$ we get 
$\M\vDash \exists z\,S(a^h,b^h,z).$ Claim proved.

Claim 2. For $A\subset N,$ $|A|< |M|,$ assume
$h_A$ exists and let $\hat b\in M\setminus h_A(A).$ Then there is a $b\in N$ and an extension $ h_{Ab}: b\mapsto \hat b.$

Proof. Consider the type over $A,$ 
$$p(y):=\{ \neg \exists z\, S(a,y,z): \M\vDash \neg \exists z\, S(a^h,\hat b,z):
 a\subseteq_\mathrm{fin} A,\ S \mbox{ positive}\}.$$
This is consistent in ${\bf N}$ since otherwise 
$${\bf N}\vDash \forall y\bigvee_{i=1}^k \exists z_i\, S_i(a,y,z_i)$$
for the negation of a finite subtype. The formula after $\vDash$
 is positive so  
$$\M\vDash \forall y\bigvee_{i=1}^k \exists z_i\, S_i(a^h,y,z_i)$$
$$\M\vDash \bigvee_{i=1}^k \exists z_i\, S_i(a^h,\hat b,z_i),$$
the contradiction. 

Now, since the type is consistent and $\mathbf{N}$ is $|M|$-saturated, there is $b\in N$ satisfying the type. Then the extension of $h_{A}$ to $h_{Ab}: b\mapsto \hat{b}$ is well defined. Claim proved. 

To finish the proof of the proposition consider a  maximal partial strong
homomorphism $h= h_A: A\to M.$  By Claim 2, $h_A$ is surjective. By Claim 1, $A=N.$ $\Box$

\epk
\bpk \label{LocGlob0} {\bf Remark.} It follows from Claim 2 in \ref{PropGlob1} that under the assumptions  any partial strong homomorphism $h_A: {^*\M}\to \M$ can be extended to a global one.
\epk
\bpk {\bf Corollary.} {\em  Assuming that $\M$ is compact, the conditions in the definition of global approximation are equivalent to the condition: 
 there exist  ${^*\M}$ and $\lm^\mathrm{glob}$ as in (\ref{M*}) and (\ref{lm2}).

Moreover,   
assuming the continuum hypothesis,
one can just take  \be \label{asM*} {^*\M}= {\prod}_I\M_i/\mathcal{D}.\ee
}

For simplicity of notation we assume (\ref{asM*}) throughout.

\epk
\bpk \label{L+} {\bf The expansion to an emerging metric language. } We are going to work in the setting of general structures as in \ref{Kies} and thus, for an arbitrary $L$-structure  
 as in \ref{global} we interpret an $n$-ary predicate $P$ on a set $M$ as a map $P^M: M^n\to \{ 0,1 \}\subset [0,1]\subset \R.$ 
 
  Let $L^+$ be the vocabulary obtained from the initial  $L$ by  adding   predicates  $\mathsf{d}_m(x_1,x_2),$ $m\in \N$ with values in $[0,1]\subset \R$  along with  unary predicates $S_m,$ $m\in \N $ (with values in $\{ 0,1 \}$).

We assume that the new predicates are interpreted in an $L^+$-structure $\M=(M; L^+)$ 
so that:

\be \label{unionM} \bigwedge_{m\in \N }  S_m(M)\subseteq S_{m+1}(M);\ee

\be \label{Mi}  \bigcup_{m\in \N }  S_m(M) =M. \ee

And for $m< k\in \N$

\be \label{Mmd} x_1,x_2\in  S_m \Rightarrow \mathsf{d}_k(x_1,x_2)=
\frac{m }{k}\mathsf{d}_m(x_1,x_2).\ee

We define a new binary predicate $\mathsf{d}(x_1,x_2)$ with values in $\R$ (possibly outside $[0,1]$ )  in terms of predicates $\mathsf{d}_m(x_1,x_2):$
$$x_1,x_2\in S_m(M) \Rightarrow \mathsf{d}(x_1,x_2):=m\mathsf{d}_m(x_1,x_2).$$

This is well-defined because of (\ref{Mmd}) and  satisfies
\be\label{mm} x_1,x_2\in S_m(M) \Rightarrow \mathsf{d}(x_1,x_2)\le m.\ee

\epk
\bpk\label{M-ultra} {\bf Local ultraproduct of emerging metric structures.} Let $\M_{i\in I}$ be a family of $L^+$-structures obtained  by expanding  $L$-structures to the emerging metric language.
 
Let $\mathcal{D}$ be a non-principal ultrafilter on $I.$ The local ultraproduct of $\M_{i\in I}$ along $\mathcal{D}$ is constructed as follows.

We treat structures $\M_i$ as two-sorted first-order structures $(M_i,\R)$ with $n$-ary predicates $P$ of vocabulary $L^+$ interpreted as maps $P: M_i^n\to [0,1]\subset \R.$

 Let ${^*\M}={\prod}_I \M_i/\mathcal{D}$  be the first-order ultraproduct of $L^+$-structures.
 
 Define ${^*\M}^\mathrm{loc},$ the {\bf local part} of ${^*\M}$ 
 to be  the substructure of ${^*\M}$ with the  $M$-universe restricted:
 $$ {^*M}^\mathrm{loc}:= {^*M}\cap \bigcup_{m\in \N}  S_m({^*M}). $$
 
We say that the local ultraproduct is well-defined if the following hold for any $x_1,\ldots, x_n,y_1,\ldots, y_n\in  {^*M}^\mathrm{loc}$:

\be \label{triangle} \mathsf{d}(x_1,x_2) + \mathsf{d}(x_2,x_3)\ge \mathsf{d}(x_1,x_3)\ee

 for each  $n$-ary $L$-predicate $P$ and  a point $(a_1...a_n)\in S_m^n({^*M}^\mathrm{loc}) \setminus P({^*M}^\mathrm{loc})$ there is $\delta\in \R_{> 0}$ such that
\be \label{deltaP}  \bigwedge_{i=1}^n \mathsf{d}(a_i,x_i)<\delta\to (x_1...x_n)\notin P \ee

In this case the relation 
$x\approx y$ defined as $$x\approx y \Leftrightarrow \bigwedge_{n\in \N }\mathsf{d}(x,y)\le \frac{1}{n}$$
is a congruence relation on  ${^*M}^\mathrm{loc}.$

The  {\bf local   ultraproduct}  of    $\M_{i\in I}$ is the $L^+$-structure $\M^\mathrm{loc}$ obtained as the quotient $$\M^\mathrm{loc}:={^*\M}^\mathrm{loc}/\approx$$ 
under the assumptions.

When the assumptions are satisfied, we say $\M_{i\in I}$  is a family with {\bf  emerging  metric}.

Note that the construction is the same as the ultraproduct construction for Keisler's general structures, \cite{KeislerGL}.


 \epk 
 \bpk \label{Los}{\bf Theorem.} {\em 
 $\M^\mathrm{loc}$ 
  is a complete metric structure with $L$-predicates defining closed subsets 
  on the Cartesian powers of the metric space.
 
 There is an $L$-homomorphism
 \be\label{lm3}\lm^\mathrm{loc}: {^*\M}^\mathrm{loc} \twoheadrightarrow \M^\mathrm{loc}\ee
onto a metric structure $\M^\mathrm{loc}.$

For a fixed $m\in \N,$ let ${^*a}$ be an $n$-tuple  in $S_m({^*\M}^\mathrm{loc})$ represented by $\{ a(i)\in M_i^n: i\in I\}.$ Let $a=\lm^\mathrm{loc}({^*a})$ and    $\sigma(x)$ be a  general structure formula  in vocabulary $L^+$ and an $n$-tuple $x$ of variables and
 quantifiers restricted to $S_m.$ 
Then the logical value of the formula in $\R$ satisfies the Lo\'s property
$$\lm^\mathrm{loc}\sigma^{\M_i}(a(i))=\sigma^{\M^\mathrm{loc}}(a).$$  
}

{\bf Proof.} (\ref{lm3}) is a direct consequence of the definition and the Lo\'s property is
by \cite{KeislerGL}, 2.2.2. $\Box$
\epk

\bpk 
Accordingly, in analogy with (\ref{lm}), we write 
\be\label{lmloc}\lm^\mathrm{loc} \M_\mathbf{i}=\M^\mathrm{loc}\ee
and call the homomorphism $\lm^\mathrm{loc}$ a {\bf local approximation} of $\M^\mathrm{loc}$ by $\M_{i\in I},$ determined by a metric and the standard part map.  
\epk
\bpk\label{remark1} {\bf Remark.} Note that ${^*\M}^\mathrm{loc}$ is a substructure of ${^*\M}$ from (\ref{lm2}). If we also assume that there is an embedding \be\label{emdloc}\M^\mathrm{loc}\hookrightarrow \M\ee then   
the map $\lm^\mathrm{loc}$ in (\ref{lm3}) is a partial homomorphism ${^*\M} \to \M.$

Under extra assumptions $\lm^\mathrm{loc}$ is a strong partial homomorphism in the sense of


\medspace

{\bf Proposition.} {\em Asssume that $\M$ is a compact $L$-topological structure.  Suppose,  for some $L$-structures $\M_{i\in I},$ (\ref{lm3}),  (\ref{emdloc}) hold and that the partial homomorphism $\lm^\mathrm{loc}$ is strong. Then the local approximation $\lm^\mathrm{loc}$ can be extended to the global approximation $\lm^\mathrm{glob}$.


}

{\bf Proof.} Immediate from \ref{LocGlob0}. \qed


\medspace

 
\epk
\bpk {\bf Approximation of metric structures.} Suppose now that $L$ is a vocabulary of an
 unbounded metric structure, see \cite{BBHU}, $L=L^+,$  with sorts $S_m$ of diameter $m,$ and $\M_i$ are 
 $L^+$-metric  structures. Then  the metric ultraproduct as defined in \cite{BBHU} is the same as  $${^*\M}^\mathrm{loc}/\approx$$
that is isomorphic to $\M^\mathrm{loc}.$

{\em  The ultraproduct of metric structures is a local approximation}.

\epk

\section{Fields and rings}
\bpk\label{3.1} We may assume, using an argument in \cite{ZCL},   that for a saturated enough model $\ZZ$ of arithmetic, there is a submodel  $ { ^\phi\Z} $ (elements of which we call {\em feasible}  numbers) and 
 non-standard positive numbers $\NN, \pp\in \ZZ,$ $\pp$ prime,  such that 
\be\label{precZ}\Z\prec  { ^\phi\Z}  \prec \ZZ\mbox{ and } \NN\ge \pp>  { ^\phi\Z}  \ee

Let $\lfr\in  { ^\phi\Z} _{>0}$ be an infinite integer. Note that by (\ref{precZ}) $\lfr^n< \pp$ for all  $n\in \N .$

Set, for $n\in \N,$ $$ { \ZZ} _{/\lfr}(n):= \{ k\in \ZZ:\  |k|<\lfr^n \},\ \  { \ZZ} _{/\lfr}:=\bigcup_{n\in \N } { ^\phi\Z} _{/\lfr}(n).$$
Note that 
$$ { \ZZ} _{/\lfr}\subseteq {^\phi\Z}.$$ 
\epk 
\bpk \label{3.2} Let $$\F=\ZZ/\pp \ZZ\cong\prod  \F_p/ _\mathcal{D}$$ be a pseudofinite field, an ultraproduct of finite prime fields $\F_p.$ Here 
  $\pp$ is the ultraproduct limit of the sequence of standard primes $p,$ that is, $\pp$ is represented in $p$-th coordinate of the ultraproduct by $p;$ $\pp(p)=p.$

By assumptions,  there is an injective map, for all $n,$
\be \label{ZfF}i_F: { \ZZ} _{/\lfr}(n)\subset \F;\ \ k\mapsto k\,\mathrm{mod}\,\pp\ee
and we will assume that it is an actual embedding.

Consider also the subsets of $\F,$ for $m\in \N ,$ 

$$S_m(\F):=\{ z\in \F: \ \exists k_1,k_2\in  { \ZZ} _{/\lfr}(m)\ 
z=k_1\cdot k_2\inv \ \&\ \frac{1}{m}\le \frac{|k_1|}{|k_2|}\le m \}$$
\be\label{sumE} \F_{/\lfr}:=\bigcup_{m\in \N }S_m(\F)\cup \{ 0\}.\ee
Note that each $S_m(\F)$ is definable in the ring $\ZZ.$

\medspace

{\bf Claim.} $\F_{/\lfr}$ is a subring of $\F.$

Indeed, $z_1,z_2\in S_m(\F)$ implies $z_1\cdot z_2 \in S_{m^2}(\F)$ and $z_1\pm z_2\in  S_{2m}(\F).$ 

\epk
\bpk
For an element $z=z(k_1,k_2)\in S_m(\F)$ as above, define 
$$||z||:= \mathrm{st}( \frac{|k_1|}{|k_2|})\in \R_{\ge 0}$$
where $|k_1|,|k_2|$ is the minimal pair such that $z=k_1\cdot k_2\inv$ and $\mathrm{st}$ is the standard part map.

This is  a $\R$-valued norm on the subring $\F_{/\lfr}$ of $\F.$ 
Define respectively $$\mathsf{d}(z_1,z_2)= || z_1-z_2||.$$
The same calculation as in the proof of the claim in \ref{3.2} shows  
that the norm and the distance are well-defined for the  pseudo-finite $\F$ with $\lfr$ as in  \ref{3.1}.

\epk
\bpk\label{3.4} {\bf Finite fields.} By definition, 
   the element $\lfr\in \ZZ$ is represented in each coordinate of the ultraproduct by $\lfr(p)\in \Z,$ the sequence $\lfr(p)$ is unique modulo the ultrafilter.

For each $n\in \N$ the embedding (\ref{ZfF}) holds for almost all $p$ with respect to $\mathcal{D}$ 
and we   
can formally apply the above definitions to $\F=\F_p,$ interpreting $\lfr$ as $\lfr(p).$ In this setting the union (\ref{sumE}) is equivalent to a finite union and $\F_{/\lfr}=\F_p.$ Note that the distance satisfies (\ref{triangle}) only for $x_1,x_2\in  S_m(\F)$ with $\lfr(p)^m$ sufficiently small compared to $p.$  

\medskip

 \label{Prop3.4}
{\bf Proposition.} {\em Finite fields $\F_p$ locally approximate the field of real numbers:
\be\label{locF} \lm ^\mathrm{loc} \F_\pp = \R.\ee 
 $\lm ^\mathrm{loc} $ can be extended to a global approximation: 
 \be\label{globF} \lm ^\mathrm{glob} \F_\pp = \bar\C\ee
where $\bar{\C}$ is the one-point compactification of the field of complex numbers. }

{\bf Proof.} The local part of an ultraproduct of finite fields  $\prod_{\pp}^{\mathrm{loc}} \F_p/\mathcal{D},$    determined by the choice of predicates $S_m$ in \ref{3.2}, is exactly  $\F_{/\lfr}.$  Factorising by the equivalence $\approx$ corresponds to application of the standard part map, which gives us $\R.$ 

The second part of the proposition follows from \ref{LocGlob0}.
$\Box$

\medspace

{\bf Remark.} In addition to (\ref{globF}), Theorem 5.2 of \cite{Perfect} proves that no other locally compact field can be globally approximated by finite fields. The compactification is also necessary, since the kernel of the homomorphism $\lm^\mathrm{glob}$ has to be non-trivial.  
\epk
\bpk {\bf Rings.}  
Consider a pseudo-finite ring $$\K=\K_\NN=\ZZ{/\NN \ZZ}\cong {\prod}  \Z_n/\mathcal{D} ,$$
the ultraproduct of finite residue rings $\Z_n,$ as topological structure in Zariski topology. 

We will assume that  $\pp|\NN$ and thus there is a surjective ring homomorphism
$$s_{K,F}: \K\to \F;\ \ k\mapsto k\,\mathrm{mod}\, \pp.$$  

On the other hand, as in (\ref{ZfF}) there is an embedding of rings
$$i_{K}: { ^\phi\Z} _{/\lfr}\hookrightarrow \K$$
and this agrees with $i_F$ of (\ref{ZfF}) in such a way that $s_{K,F}\circ i_K=i_F.$ 

Define $$K_{/\lfr} := s_{K,F}\inv \F_{/\lfr}$$
and use pseudo-metric $\mathsf{d}=\mathsf{d}_F$ on $\F_{/\lfr}$ to define a pseudo-metric on $\K_{/\lfr}$ $$\mathsf{d}_K(x,y):=\mathsf{d}_F(s_{K,F}(x), s_{K,F}(y)).$$

Set, for $m\in \N ,$ $$\K_{/\lfr}(m):= \{ z\in \K_{/\lfr}: \mathsf{d}(0,z)\le m \}.$$

With these notations we can proceed with the same construction and arguments as in the case of fields and obtain
\be\label{locK} \lm ^\mathrm{loc} \K_n = \R.\ee 
In the global setting we obtain again
\be\label{globK} \lm ^\mathrm{glob} \K_n = \bar\C\ee
due to the fact that $\lm^\mathrm{loc} $  filters through $\K\twoheadrightarrow \F,$
\epk 
\bpk\label{3.6} {\bf The Lorentz-invariant space-time.} A finite approximation to a Lorentz-invariant space-time was constructed in \cite{ZlbKL}.
 Theorem 3.7 (19) therein in the current notation states  $$ \lm ^\mathrm{loc} \left( \MM(\K_n), \mathrm{SL}(2,\K_n^{(2)})\right)  = \left( \MM(\R), \mathrm{SL}(2,\C)\right), $$
 where $\MM(\R)$ is $\R^4$ equipped with  Minkowski metric $t^2-x^2-y^2-z^2$ and $\mathrm{SL}(2,\C)$ acts on $\R^4$ preserving the metric. Similarly $\MM(\K_n)$ is $K_n^4,$ for a pseudo-finite ring as in (\ref{locK}), with a ``metric'' given by the same formula, and $$\K_n^{(2)}= \{ (a,b)\in \K^2: \begin{array}{ll} (a_1,b_1)+(a_2,b_2)=(a_1+a_2, b_1+b_2); \\ (a_1,b_1)\cdot (a_2,b_2)=(a_1\cdot a_2 - b_1\cdot b_2,\, a_1\cdot b_2 + b_1\cdot a_2)\end{array}   \}.  $$  

Recalling that the 2-element centre $C_2$ of $\mathrm{SL}(2,\C)$ acts on $\MM(\R)$ trivially and that 
$\mathrm{SL}(2,\C)/C_2\cong \mathrm{SO}^+(1,3),$ we establish the equivalent statement in terms of the Lorentz group $\mathrm{SO}^+(1,3):$
$$ \lm ^\mathrm{loc} \left( \MM(\K_n), \mathrm{SL}(2,\K_n^{(2)})/C_2\right)  = \left( \MM(\R), \mathrm{SO}^+(1,3)\right). $$
That is Minkowski space equipped with Minkowski metric  and the invariant action of Lorentz group on it is locally approximable by finite   spaces $\MM(\K_n)$ equipped with an analogue of Minkowski metric  with the invariant action of finite groups  $\mathrm{SL}(2,\K_n^{(2)})/C_2.$  
\epk
\bpk {\bf The ring $\Z$ is locally approximated by finite rings.}

Indeed, restrict the definition of $S_m$ to $$S_m^\Z(\F)=\{ z\in \F: \ z=k\in  { ^\phi\Z} _{/\lfr}(m), \ |k|\le m \} \subset S_m(\F) $$ and define the distance in the same way as in \ref{3.2}.  Then the local structure defined by  $S_m^\Z$ is a substructure of $\F^\mathrm{loc}$ and, for  the same $\lm^\mathrm{loc},$ $$\lm^\mathrm{loc}: \bigcup_m S_m^\Z(\F) \twoheadrightarrow \Z.$$ $\Box$ 
\epk
\bpk\label{QVar} {\bf Varieties over fields and rings.} Let $V\subseteq \mathbb{A}^n$ be an affine variety defined over $\Z,$
$V(\F_p)$ its set of $\F_p$-points. We want to calculate $\lm^\mathrm{loc} V(\F_p),$ 
where $\lm^\mathrm{loc}$ on $\F_p$-points is the same as in (\ref{locF}). In other words, the structure in question is an affine variety over a field and the emerging metric is determined on the fields in \ref{3.2} - \ref{3.4}.   

{\bf Proposition.} $$\lm^\mathrm{loc} V(\F_p)=\overline{V(\Q)}\subseteq V(\R)$$ 
{\em $\overline{V(\Q)}$ is the set of  
	limit points of the set  of $\Q$-ponits in $V(\R).$ }

{\bf Proof.} For the pseudofinite $\F_\pp=\F$ we need to consider the image of $V(\F_{/\lfr}).$ The argument in the proof of \ref{Prop3.4} tells us that this is the same as the image of $V(\QQ_\mathrm{fin})$ under the standard part map, which gives us $ \overline{V(\Q)}.$ $\Box$

\epk
\bpk {\bf Corollary.} With the same notation, for finite rings $\K_n$ as above, $$\lm^\mathrm{loc} V(\K_n)=\overline{V(\Q)}\subseteq V(\R).$$
\epk
\bpk  The abelian multiplicative group $\Ss^1=\{ q\in \C^\times: \ |q|=1\}\cong U(1)$  is both locally and globally approximated by finite cyclic groups. For the local approximation one uses the Lie group metric on $U(1)$ and its finite cyclic subgroups. Since the  topology defined by the metric is compact, the same homomorphism $\lm^\mathrm{loc}$ serves as a global approximation $\lm^\mathrm{glob}.$ 

\epk

\bpk \label{Heis}
The   Heisenberg group $\mathrm{H}(3)$ is the connected, simply  connected Lie group with the Lie algebra generated  by self-adjoint operators (on a Hilbert space) $\Pp$ and $\Qq$ satisfying the identity 
$$\Pp\Qq- \Qq\Pp=i\hh\mathrm{I}$$
for some generic $\hh\in \R.$
Representations of this group are the subject of the Stone-von Neumann theorem. Every element of the group can be written in the form $$ q\cdot \e^{ai \Pp}\cdot \e^{bi \Qq}, \mbox{ with } \ a,b\in \R,\  q\in \Ss^1 .$$

It is easy to see that for $\hh=r\in 2\pi \Q$ one gets a finite group acting on a finite-dimensional space. Such finite groups  locally  approximate
 $\mathrm{H}(3)$ in the metric of $\Ss^1\times \R\times \R.$

\epk
\section{Compact simple Lie groups} \label{s4}
\bpk\label{some} We first consider some well known cases.

A. $\mathrm{SO}_n(\R),$ $n>2.$ It is well known, see e.g., \cite{Borel}, 18.2, that the rational points form the dense subgroup   $\mathrm{SO}_n(\Q)$ in each such group. Using  \ref{QVar} we get local approximation by finite algebraic groups
$$\lm^\mathrm{loc} \mathrm{SO}_n(\F_\pp)=  \mathrm{SO}_n(\R).$$ 

B. $\mathrm{SU}_n(\C),$ $n\ge 2.$ Again, it is easy to see that  
$\mathrm{SU}_n(\Q[i])$ is dense in  $\mathrm{SU}_n(\C).$ 

The general theory of density of $\kk$-points, which covers all compact simple Lie groups, is presented in \cite{PR}.

{\bf General Fact.} (\cite{PR}, Chapter 6) {\em Let $\kk$ be a number field, $v$ an archimedean place, and $G$ an absolutely simple algebraic group over $\kk$.
 Assume $G(\kk_v)$ is compact (equivalently: $G$ is anisotropic over $\kk_v$).
	Then $G(\kk)$ is dense in $G(\kk_v)$ in the real/complex topology.}

\medspace 

Consider the rings $\F_p^{(2)}=\F_p\times \F_p$ with operations 
$$\begin{array}{ll} (a_1,b_1)+(a_2,b_2)= (a_1+a_2, b_1+b_2),\\   (a_1,b_1)\cdot(a_2,b_2)= (a_1\cdot a_2- b_1\cdot b_2, a_1\cdot b_2+b_1\cdot a_2)\end{array}$$ 
It is clear that, using the same emerging metric on $\F_p$, 
$$\lm^\mathrm{loc} \F_\pp^{(2)}= \R[i]=  \C,$$
and so
$$\lm^\mathrm{loc} \mathrm{SU}_n(\F_\pp^{(2)})=  \mathrm{SU}_n(\C).$$


\epk

\bpk
{\bf Theorem.} {\em Every simple compact Lie group is locally approximated by finite groups.}

{\bf Proof.} It is enough to establish that each compact simple group $G$ satisfies the $\kk$-rationality condition of \ref{QVar} for $\kk=\Q$ or $\kk=\Q[i].$ The classification of such groups is well known, see \cite{SLG}:

1. $\mathrm{SO}(n)$, $n\ge 3$ (types $B_n, D_n$) 

$\kk=\Q,$ the archimedean place $v = \infty$, so $\kk_v = \R.$

This is the case A of \ref{some}.

\medspace

2. $\mathrm{SU}(n+1)$ (compact type $A_n$)

The compact  $\mathrm{SU}(n+1)$  is most naturally realised as a unitary group over an imaginary quadratic field. This is the case B of \ref{some} above.

\medspace

3. $\mathrm{Sp}(n)$ (type $C_n$, the compact symplectic group)

Here one must be a bit more careful: the split symplectic group $\mathrm{Sp}_{2n}$ 
over $\Q$ has non‑compact real points  $\mathrm{Sp}_{2n}(\R).$ The compact 
$\mathrm{Sp}(n)$ is an inner form, i.e. is representable as the group of isometries of a nondegenerate Hermitian form over the quaternion division algebra of rank $n.$

The field is $\kk = \Q$, and $G(\Q)$ is dense in   $\mathrm{Sp}(n).$

\medspace

4. Exceptional compact groups

For the compact real forms of $G_2, F_4, E_6, E_7, E_8$, there are anisotropic $\Q$-forms whose real points are exactly the compact groups listed above. One convenient choice in each case is:

Field $\kk = \Q.$

Group: An absolutely simple, simply connected $\Q$-group of the given type, anisotropic over $\R.$ 

Concretely:

$G_2$: Automorphism group of a rational octonion division algebra (a Cayley algebra) over $\Q.$ Then $G(\R)$ is the compact $G_2.$

$F_4$: Automorphism group of an Albert division algebra over $\Q$, anisotropic over $\R$;

$E_6, E_7, E_8$:  By Tits's classification of $\Q$-forms there exist simply connected $\Q$-groups of each type that are anisotropic at $\infty$; their real points are the compact real forms $G(\R).$

$\Box$

\medspace

The theorem  is  in contrast with the fact proved in \cite{Nik} that the global approximation of simple Lie groups by finite groups is not possible.

\epk

\end{document}